\documentstyle[twoside,amstex]{article}
\input amssym.def
\input amssym.tex
\def\bc{\begin{center}}
\def\ec{\end{center}}
\def\no{\noindent}

\begin{document}
\thispagestyle{empty} \vspace*{3 true cm} \pagestyle{myheadings}
\markboth {\hfill {\sl Huanyin Chen, S. Halicioglu and H.
Kose}\hfill} {\hfill{\sl On Perfectly Clean Rings}\hfill}
\vspace*{-1.5 true cm} \bc{\large\bf On Perfectly Clean Rings}\ec

\vskip6mm
\bc{{\bf Huanyin Chen}\\[1mm]
Department of Mathematics, Hangzhou Normal University\\
Hangzhou 310036, China, huanyinchen@@aliyun.com}\ec
\bc{{\bf S. Halicioglu}\\[1mm]
Department of Mathematics, Ankara University\\ 06100 Ankara,
Turkey, halici@@science.ankara.edu.tr}\ec
\bc{{\bf H. Kose}\\[1mm]
Department of Mathematics, Ahi Evran University\\ Kirsehir,
Turkey, handankose@@gmail.com}\ec

\begin{abstract} An element $a$
of a ring $R$ is called perfectly clean if there exists an
idempotent $e\in comm^2(a)$ such that $a-e\in U(R)$. A ring $R$ is
perfectly clean in case every element in $R$ is perfectly clean.
In this paper, we investigate conditions on a local ring $R$ that
imply that $2\times 2$ matrix rings and triangular matrix rings
are perfectly clean. We shall show that for these rings perfect
cleanness and strong cleanness coincide with each other, and
enhance many known results. We also obtain several criteria for
such a triangular matrix ring to be perfectly $J$-clean. For
instance, it is proved that for a commutative ring $R$, $T_n(R)$
is perfectly $J$-clean if and only if $R$ is strongly $J$-clean.\\[2mm]
\end{abstract}

\section{Introduction}

The commutant and double commutant of an element $a$ in a ring $R$
are defined by $comm(a)=\{x\in R~|~xa=ax\}$, $comm^2(a)=\{x\in
R~|~xy=yx~\mbox{for all}~y\in comm(a)\}$, respectively. An element
$a\in R$ is strongly clean provided that there exists an
idempotent $e\in comm(a)$ such that $a-e\in U(R)$. A ring $R$ is
called strongly clean in the case that every element in $R$ is
strongly clean. Strongly clean matrix rings and triangular matrix
rings over local rings are extensively studied by many authors
(cf. [1-2], [5-6] and [12-13]. An element $a\in R$ is quasipolar
provided that there exists an idempotent $e\in comm^2(a)$ such
that $a+e\in U(R)$ and $ae\in R^{qnil}$, while $R^{qnil}=\{ x\in
R~|~1+xr\in U(R)~\mbox{for any}~r\in comm(x)\}$. A ring $R$ is
called quasipolar if every element in $R$ is quasipolar. As is
well known, a ring $R$ is quasipolar if and only if for any $a\in
R$ there exists a $b\in comm^2(a)$ such that $b=bab$ and
$b-b^2a\in R^{qnil}$. This concept evolved from Banach algebra. In
fact, for a Banach algebra $R$, $$a\in R^{qnil}\Leftrightarrow
\lim\limits_{n\to\infty}\parallel a^n\parallel^{\frac{1}{n}}=0.$$
It is shown that every quasipolar ring is strongly clean.
Recently, quasipolar $2\times 2$ matrix rings and triangular
matrix rings over local rings are also studied from different
point of views (cf. [7-9] and [11]).

The motivation of this article is to introduce a medium class
between strongly clean rings and quasipolar rings. An element $a$
of a ring $R$ is called perfectly clean if there exists an
idempotent $e\in comm^2(a)$ such that $a-e\in U(R)$. A ring $R$ is
perfectly clean in case every element in $R$ is perfectly clean.
We shall show that for $2\times 2$ matrix rings and triangular
matrix rings over local rings, perfect cleanness and strong
cleanness coincide with each other, and enhance many known
results, e.g. [5, Theorem 8], [11, Theorem 2.8] and [12, Theorem
7]. Replaced $U(R)$ by $J(R)$, we introduce perfect $J$-clean
rings as a subclass of perfectly clean rings. Furthermore, we show
that strong $J$-cleanness for triangular matrix rings over a local
ring can be enhanced to such stronger properties. These also
generalize the corresponding properties of $J$-quasipolarity,
e.g., [8, Theorem 4.9].

We write $U(R)$ and $J(R)$ for the set of all invertible elements
and the Jacobson radical of $R$. $M_n(R)$ and $T_n(R)$ stand for
the rings of all $n\times n$ matrices and triangular matrices over
a ring $R$.

\section{Elementary Properties}

\vskip4mm Clearly, abelian exchange ring is perfectly clean. Every
quasipolar ring is perfectly clean. For instance, every strongly
$\pi$-regular ring. In fact, we have $\{$ quasipolar rings
$\}\subsetneq \{$ perfectly clean rings $\}\subsetneq \{$ strongly
clean rings $\}$. We begin with

\vskip4mm \hspace{-1.8em} {\bf Theorem 2.1.}\ \ {\it Let $R$ be a
ring. Then the following are equivalent:} \vspace{-.5mm}
\begin{enumerate}
\item [(1)]{\it $R$ is perfectly clean.}
\vspace{-.5mm}
\item [(2)]{\it For any $a\in R$, there exists an $x\in comm^2(a)$ such that $x=xax$ and
$1-x\in (1-a)R\bigcap R(1-a)$.}
\end{enumerate}
\vspace{-.5mm} {\it Proof.}\ \ $(1)\Rightarrow (2)$ For any $a\in
R$, there exists an idempotent $e\in comm^2(a)$ such that
$u:=a-e\in U(R)$. Set $x=u^{-1}(1-e)$. Let $y\in comm(a)$. Then
$ay=ya$. As $uy=(a-e)y=y(a-e)=yu$, we get $u^{-1}y=yu^{-1}$. Thus,
$xy=u^{-1}(1-e)y=u^{-1}y(1-e)=yu^{-1}(1-e)=yx$. This implies that
$x\in comm^2(a)$. Further,
$xax=u^{-1}(1-e)(u+e)u^{-1}(1-e)=u^{-1}(1-e)=x$. Clearly,
$u=(1-e)-(1-a)$, and so $1-u^{-1}(1-e)=u^{-1}(1-a)$. This implies
that $1-x\in R(1-a)$. Likewise, $1-x\in (1-a)R$ as
$(1-e)u^{-1}=u^{-1}(1-e)$. Therefore $1-x\in (1-a)R\bigcap
R(1-a)$, as required.

$(2)\Rightarrow (1)$ For any $a\in R$, there exists an $x\in
comm^2(a)$ such that $x=xax$ and $1-x\in (1-a)R\bigcap R(1-a)$.
Write $e=1-ax$. If $y\in comm(a)$, then $ay=ya$, and so
$axy=ayx=yax$. This shows that $ey=ye$; hence, $e\in comm^2(a)$.
In addition, $ex=xe=0$. Write $1-x=(1-a)s=t(1-a)$ for some $s,t\in
R$. Then
$$\begin{array}{lll}
(a-e)(x-es)&=&ax-aes+es\\
&=&ax+(1-a)es\\
&=&ax+e(1-a)s\\
&=&ax+e(1-x)\\
&=&ax+e\\
&=&1. \end{array}$$ Likewise, $(x-te)(a-e)=1$. Therefore $a-e\in
U(R)$, as desired.\hfill$\Box$

\vskip4mm \hspace{-1.8em} {\bf Corollary 2.2.}\ \ {\it Let $R$ be
a ring. Then the following are equivalent:}\vspace{-.5mm}
\begin{enumerate}
\item [(1)]{\it $R$ is perfectly clean.}
\vspace{-.5mm}
\item [(2)]{\it For any $a\in R$, there exists an
idempotent $e\in comm^2(a)$ such that $eae\in U(eRe)$ and
$(1-e)(1-a)(1-e)\in U\big((1-e)R(1-e)\big)$.}\vspace{-.5mm}
\end{enumerate} {\it Proof.}\ \ $(1)\Rightarrow (2)$ For any $a\in
R$, it follows from Theorem 2.1 that there exists an $x\in
comm^2(a)$ such that $x=xax$ and $1-x\in (1-a)R\bigcap R(1-a)$.
Write $1-x=(1-a)s=t(1-a)$ for some $s,t\in R$. Set $e=ax$. For any
$y\in comm(a)$, we have $ay=ya$, and so
$ey=(ax)y=a(yx)=(ay)x=y(ax)=ye$. Hence, $e^2=e\in comm^2(a)$.
Clearly, $(eae)(exe)=(exe)(eae)=e$; hence, $eae\in U(eRe)$.
Furthermore, $1-e=(1-x)+(1-a)x=(1-a)(s+x)$. This shows that
$(1-e)(1-a)(1-e)(1-x)(1-e)=1-e$. Likewise,
$(1-e)(1-x)(1-e)(1-e)(1-a)(1-e)=1-e$. Therefore
$(1-e)(1-a)(1-e)\in U\big((1-e)R(1-e)\big)$.

$(2)\Rightarrow (1)$ For any $a\in R$, we have an idempotent $e\in
comm^2(a)$ such that $eae\in U(eRe)$ and $(1-e)(1-a)(1-e)\in
U\big((1-e)R(1-e)\big)$. Hence,
$a-(1-e)=\big(eae-(1-e)(1-a)(1-e)\big)\in U(R)$. Set $p=1-e$. Then
$a-p\in U(R)$ with $p\in comm^2(a)$, as desired.\hfill$\Box$

\vskip4mm Recall that a ring $R$ is strongly nil clean provide
that every element in $R$ is the sum of an idempotent and a
nilpotent element that commutate (cf. [4] and [10]).

\vskip4mm \hspace{-1.8em} {\bf Theorem 2.3.}\ \ {\it Let $R$ be a
ring. Then $R$ is strongly nil clean if and only if}
\vspace{-.5mm}
\begin{enumerate}
\item [(1)]{\it $R$ is perfectly clean;}
\vspace{-.5mm}
\item [(2)]{\it $N(R)=\{ x\in R~|~1-x\in U(R)\}$.}
\end{enumerate}
\vspace{-.5mm} {\it Proof.}\ \ Suppose $R$ is strongly nil clean.
For any $a\in R$, we see that $a-a^2\in N(R)$. Write
$(a-a^2)^n=0$. Let $f(t)=\sum\limits_{i=0}^n \left(
\begin{array}{c}
2n\\
i
\end{array}
\right)t^{2n-i}(1-t)^i\in {\Bbb Z}[t]$. Then $f(t)\equiv 0
\,(mod\,t^n)$. Clearly, $$f(t)+ \sum\limits_{i=n+1}^{2n}
\left( \begin{array}{c} 2n\\
i \end{array} \right)x^{2n-i}(1-t)^i=\big(t+(1-t)\big)^n=1;$$
hence, $f(t)\equiv 1 \,\big(mod\,(1-t)^n\big)$. This shows that
$f(t)\big(1-f(t)\big)\equiv 0\,\big(mod\,t^n(1-t)^n\big)$. Let
$e=f(a)$. Then $e\in R$ is an idempotent. For any $x\in comm(a)$,
we see that $xa=ax$, and so $xe=xf(a)=f(a)x=ex$. Thus, $x\in
comm^2(a)$. Furthermore,
$a-e=a-a^{2n}+(a-a^2)g(a)=(a-a^2)\big(1+a+a^2+\cdots
+a^{2n-2}+g(a)\big)\in N(R)$, where $g(t)\in {\Bbb Z}[t]$. Thus,
$a=(1-e)+(2e-1+a-e)$ with $1-e\in comm^2(a)$ and $2e-1+a-e\in
U(R)$. Therefore, $R$ is perfectly clean. Clearly, $N(R)\subseteq
\{ x\in R~|~1-x\in U(R)\}$. If $1-x\in U(R)$, then $x=e+w$ with
$e\in comm(x)$ and $w\in W(R)$. Hence, $1-e=(1-x)+w\in U(R)$. This
implies that $1-e=1$, and so $x=w\in N(R)$. Therefore $N(R)=\{
x\in R~|~1-x\in U(R)\}$.

Conversely, assume that $(1)$ and $(2)$ hold. For any $a\in R$,
there exists an idempotent $e\in comm^2(a)$ and a unit $u\in R$
such that $-a=e-u$. Hence, $a=-e+u=(1-e)-(1-u)$. By hypothesis,
$1-u\in N(R)$. Accordingly, $R$ is strongly nil clean.\hfill$\Box$

\vskip4mm \hspace{-1.8em} {\bf Corollary 2.4.}\ \ {\it Let $R$ be
a ring. Then $R$ is strongly nil clean if and only if}
\vspace{-.5mm}
\begin{enumerate}
\item [(1)]{\it $R$ is quasipolar;}
\vspace{-.5mm}
\item [(2)]{\it $N(R)=\{ x\in R~|~1-x\in U(R)\}$.}
\end{enumerate}
\vspace{-.5mm} {\it Proof.}\ \ Suppose that $R$ is strongly nil
clean. Then $(2)$ holds from Theorem 2.5. For any $a\in R$, as in
the proof of Theorem 2.5, $a=e+w$ with $e\in comm^2(a)$ and $w\in
N(R)$. Hence, $a=(1-e)+(2e-1+w)$ where $2e-1+w\in U(R)$.
Furthermore, $(1-e)a=(1-e)w\in N(R)\subseteq R^{qnil}$. Therefore
$R$ is quasipolar.

Conversely, assume that $(1)$ and $(2)$ hold. Then $R$ is
perfectly clean. Accordingly, we complete the proof by Theorem
2.3. \hfill$\Box$

\vskip4mm \hspace{-1.8em} {\bf Lemma 2.5.}\ \ {\it Let $R$ be a
ring. Then the following are equivalent:\/} \vspace{-.5mm}
\begin{enumerate}
\item [(1)]{\it $R$ is perfectly clean.}
\vspace{-.5mm}
\item [(2)]{\it For each $a\in R$ there exists an idempotent $e\in comm^2(a)$ such that $a-e$ and $a+e$ are invertible.}
\end{enumerate}
\vspace{-.5mm}  {\it Proof.}\ \ $(1)\Rightarrow (2)$ Let $a\in R$.
Then $a^2\in R$ is perfectly clean. Thus, we can find an
idempotent $e\in comm^2(a^2)$ such that $a^2-e\in U(R)$. As $a\in
comm(a^2)$, we see that $ae=ea$. Hence, $a^2-e=(a-e)(a+e)$, and
therefore we conclude that $a-e,a+e\in U(R)$.

$(2)\Rightarrow (1)$ is trivial.\hfill$\Box$

\vskip4mm \hspace{-1.8em} {\bf Theorem 2.6.}\ \ {\it Let $R$ be
perfectly clean. Then for any $A\in M_n(R)$ there exist $U,V\in
GL_n(R)$ such that $2A=U+V$.}\vskip2mm \hspace{-1.5em}{\it
Proof.}\ \ We prove the result by induction on $n$. For any $a\in
R$, there exists an idempotent $e\in comm^2(a)$ such that
$u:=a-e,v:=a+e\in U(R)$, by Lemma 6.4.24. Hence, $2a=u+v$, and so
the result holds for $n=1$. Assume that the result holds for
$n\leq k$ $(k\geq 1)$. Let $n=k+1$, and let $A\in M_n(R)$. Write
$A=\left(
\begin{array}{cc}
x&\alpha\\
\beta&X
\end{array}
\right)$, where $x\in R, \alpha\in M_{1\times k}(R), \beta\in
M_{k\times 1}(R)$ and $X\in M_{k}(R)$. In view of Lemma 6.4.24, we
have a $u\in U(R)$ such that $2x-u=v\in U(R)$. By hypothesis, we
have a $U\in GL_k(R)$ such that $2\big(X-2\beta
v^{-1}\alpha\big)-U=V\in GL_k(R)$. Hence
$$
2A-\left(
\begin{array}{cc}
u&0\\
0&U
\end{array}
\right) = \left(
\begin{array}{cc}
v&2\alpha\\
2\beta&V+4\beta v^{-1}\alpha
\end{array}
\right).$$ It is easy to verify that
$$\left(
\begin{array}{cc}
v&2\alpha\\
2\beta&V+4\beta v^{-1}\alpha
\end{array}
\right)=\left(
\begin{array}{cc}
1&\\
2\beta v^{-1}&I_k
\end{array}
\right)\left(
\begin{array}{cc}
v&2\alpha\\
0&V
\end{array}
\right)\in GL_n(R).$$ By induction, we complete the
proof.\hfill$\Box$

\vskip4mm \hspace{-1.8em} {\bf Corollary 2.7.}\ \ {\it Let $R$ be
a quasipolar ring. If $\frac{1}{2}\in R$, then every $n\times n$
matrix over $R$ is the sum of two invertible matrices.}\vskip2mm
\hspace{-1.5em}{\it Proof.}\ \ As every quasipolar ring is
perfectly, we complete the proof by Theorem 2.6.\hfill$\Box$

\vskip4mm As a consequence, we derive the following known fact:
Let $R$ be a strongly $\pi$-regular ring with $\frac{1}{2}\in R$.
Then every $n\times n$ matrix over $R$ is the sum of two
invertible matrices.

\section{Matrix Rings}

\vskip4mm Recall that a ring $R$ is local if it has only one
maximal right ideal. A ring $R$ is local if and only if for any
$a\in R$ either $a$ or $1-a$ is invertible. The necessary and
sufficient conditions under which $2\times 2$ matrix ring over a
local ring are attractive. The purpose of this section is to
enhance the known results to perfect cleanness under the same
conditions.

\vskip4mm \hspace{-1.8em} {\bf Lemma 3.1.}\ \ {\it Let $R$ be a
ring, and let $u\in U(R)$. Then the following are equivalent:}
\vspace{-.5mm}
\begin{enumerate}
\item [(1)]{\it $a\in R$ is perfectly clean.}
\vspace{-.5mm}
\item [(2)]{\it $uau^{-1}\in R$ is perfectly clean.}
\end{enumerate}
\vspace{-.5mm} {\it Proof.}\ \ $(1)\Rightarrow (2)$ By hypothesis,
there exists an idempotent $e\in comm^2(a)$ such that $a-e\in
U(R)$. Hence, $uau^{-1}-ueu^{-1}\in U(R)$. For any $x\in
comm(uau^{-1})$, we see that
$x\big(uau^{-1}\big)=\big(uau^{-1}\big)x$, and so
$$\big(u^{-1}xu\big)a=a\big(u^{-1}xu\big).$$ Thus,
$$\big(u^{-1}xu\big)e=e\big(u^{-1}xu\big).$$ Hence,
$$x\big(ueu^{-1}\big)=\big(ueu^{-1}\big)x.$$ We conclude that
$ueu^{-1}\in comm^2(uau^{-1})$, as required.

$(1)\Rightarrow (2)$ is symmetric.\hfill$\Box$

\vskip4mm A ring is weakly cobleached provided that for any $a\in
J(R),b\in 1+J(R)$, $l_a-r_b$ and $l_b-r_a$ are both injective. For
instance, every commutative local ring, every local ring with nil
Jacobson radical.

\vskip4mm \hspace{-1.8em} {\bf Theorem 3.2.}\ \ {\it Let $R$ be a
weakly cobleached local ring. Then the following are equivalent:}
\vspace{-.5mm}
\begin{enumerate}
\item [(1)]{\it $M_2(R)$ is perfectly clean.}
\vspace{-.5mm}
\item [(2)]{\it $M_2(R)$ is strongly clean.}
\vspace{-.5mm}
\item [(3)]{\it For any $A\in M_2(R)$, $A\in GL_2(R)$, or $I_2-A\in GL_2(R)$, or $A$ is
similar to a diagonal matrix.}
\end{enumerate}
\vspace{-.5mm} {\it Proof.}\ \ $(1)\Rightarrow (2)$ is trivial.

$(2)\Rightarrow (3)$ is obtained by [13, Theorem 7].

$(3)\Rightarrow (1)$ For any $A\in M_2(R)$, $A\in GL_2(R)$, or
$I_2-A\in GL_2(R)$, or $A$ is similar to a diagonal matrix. If $A$
or $I_2-A\in GL_2(R)$, then $A$ is perfectly clean. Assume now
that $A$ is similar to a diagonal matrix with $A,I_2-A\not\in
G_2(R)$. We may assume that $A$ is similar to $\left(
\begin{array}{cc}
\lambda&0\\
0&\mu \end{array} \right)$, where $\lambda\in U(R), \mu\in J(R)$.
If $\lambda\in 1+U(R)$, then $\left(
\begin{array}{cc}
\lambda&0\\
0&\mu \end{array} \right)-I_2\in GL_2(R)$; hence, it is perfectly
clean. In view of Lemma 3.1, $A$ is perfectly clean. Thus, we
assume that $\lambda\in 1+J(R)$. By Lemma 3.1, it will suffice to
show that $\left(
\begin{array}{cc}
\lambda&0\\
0&\mu \end{array} \right)\in GL_2(R)$ is perfectly clean. Clearly,
$$\left(
\begin{array}{cc}
\lambda&0\\
0&\mu \end{array} \right)=\left(
\begin{array}{cc}
0&0\\
0&1\end{array} \right)+\left(
\begin{array}{cc}
\lambda&0\\
0&\mu-1\end{array} \right),$$ where $\left(
\begin{array}{cc}
\lambda&0\\
0&\mu-1\end{array} \right)\in GL_2(R)$. We shall show that the
idempotent $\left(
\begin{array}{cc}
0&0\\
0&1\end{array} \right)\in comm^2\big(\left(
\begin{array}{cc}
\lambda&0\\
0&\mu \end{array} \right)\big)$. For any $\left(
\begin{array}{cc}
x&s\\
t&y \end{array} \right)\in comm\big(\left(
\begin{array}{cc}
\lambda&0\\
0&\mu \end{array} \right)\big)$, we see that $$\lambda
s=s\mu~\mbox{and}~\mu t=t\lambda;$$ hence, $s=t=0$. This implies
that $$\left(
\begin{array}{cc}
x&s\\
t&y \end{array} \right)\left(
\begin{array}{cc}
0&0\\
0&1\end{array} \right)=\left(
\begin{array}{cc}
0&0\\
0&y\end{array} \right)=\left(
\begin{array}{cc}
0&0\\
0&1\end{array} \right)\left(
\begin{array}{cc}
x&s\\
t&y \end{array} \right).$$ Therefore $\left(
\begin{array}{cc}
0&0\\
0&1\end{array} \right)\in comm^2\big(\left(
\begin{array}{cc}
\lambda&0\\
0&\mu \end{array} \right)\big)$, hence the result. \hfill$\Box$

\vskip4mm \hspace{-1.8em} {\bf Corollary 3.3.}\ \ {\it Let $R$ be
a commutative local ring. Then the following are equivalent:}
\vspace{-.5mm}
\begin{enumerate}
\item [(1)]{\it $M_2(R)$ is perfectly clean.}
\vspace{-.5mm}
\item [(2)]{\it $M_2(R)$ is strongly clean.}
\vspace{-.5mm}
\item [(3)]{\it For any $A\in M_2(R)$, $A\in GL_2(R)$, or $I_2-A\in GL_2(R)$, or $A$ is
similar to a diagonal matrix.}
\end{enumerate}
\vspace{-.5mm} {\it Proof.}\ \ It is immediate from Theorem 3.1 as
every commutative local ring is weakly cobleached.\hfill$\Box$

\vskip4mm Let $p$ be a prime. We use $\widehat{{\Bbb Z}_{p}}$ to
denote the ring of all $p$-adic integers. In view of [6, Theorem
2.4], $M_2\big(\widehat{{\Bbb Z}_{p}}\big)$ is strongly clean, and
therefore $M_2\big(\widehat{{\Bbb Z}_{p}}\big)$ is perfectly
clean, by Corollary 3.3.

\vskip4mm \hspace{-1.8em} {\bf Theorem 3.4.}\ \ {\it Let $R$ and
$S$ be local rings. Then the following are equivalent:}
\vspace{-.5mm}
\begin{enumerate}
\item [(1)]{\it $\left(
\begin{array}{cc}
R&V\\
0&S \end{array} \right)$ is perfectly clean.} \vspace{-.5mm}
\item [(2)]{\it For any $a\in J(R),b\in 1+J(S), v\in V$,
there exists a unique $x\in V$ such that $ax-xb=v$.}
\end{enumerate}
\vspace{-.5mm} {\it Proof.}\ \ $(1)\Rightarrow (2)$ Let $a\in
1+J(R), b\in J(R)$. Set $A=\left(
\begin{array}{cc}
a&-v\\
0&b \end{array} \right)$. By hypothesis, we can find an idempotent
$E\in comm^2(A)$ such that $A-E\in \left(
\begin{array}{cc}
R&V\\
0&S \end{array} \right)$ is invertible. Clearly, $E=\left(
\begin{array}{cc}
0&x\\
0&1 \end{array} \right)$ for some $x\in V$. Thus, $ax-xb=v$.
Suppose that $ay-yb=v$ for a $y\in V$. Then
$$A\left(
\begin{array}{cc}
0&y\\
0&1 \end{array} \right)=\left(
\begin{array}{cc}
0&y\\
0&1 \end{array} \right)A,$$ and so $\left(
\begin{array}{cc}
0&y\\
0&1 \end{array} \right)\in comm(A)$. This implies that $$E\left(
\begin{array}{cc}
0&y\\
0&1 \end{array} \right)=\left(
\begin{array}{cc}
0&y\\
0&1 \end{array} \right)E;$$ hence, $x=y$. Therefore there exists a
unique $x\in V$ such that $ax-xb=v$, as desired.

$(2)\Rightarrow (1)$ Let $T=\left(
\begin{array}{cc}
R&V\\
0&S \end{array} \right)$, and let $A=\left(
\begin{array}{cc}
a&v\\
0&b \end{array} \right)\in \left(
\begin{array}{cc}
R&V\\
0&S \end{array} \right)$.

Case I. $a\in J(R),b\in J(S)$. Then $A-\left(
\begin{array}{cc}
1_R&0\\
0&1_S \end{array} \right)\in U(T)$; hence, $A$ is perfectly clean.

Case II. $a\in U(R),b\in U(S)$. Then $A-0\in U(T)$; hence, $A$ is
perfectly clean.

Case III. $a\in U(R), b\in J(S)$. $(i)$ $a\in 1+U(R), b\in J(S)$.
Then $A-\left(
\begin{array}{cc}
1_R&0\\
0&1_S \end{array} \right)\in T$ is invertible; hence, $A\in T$ is
perfectly. $(ii)$ $a\in 1+J(R), b\in J(S)$. Then we can find a
$t\in V$ such that $at-tb=-v$. Let $\left(
\begin{array}{cc}
x&s\\
0&y \end{array} \right)\in comm(A)$. Then
$$A\left(
\begin{array}{cc}
x&s\\
0&y \end{array} \right)=\left(
\begin{array}{cc}
x&s\\
0&y \end{array} \right)A,$$ and so $$ax=xa,by=yb, as-sb=xv-vy.$$
Hence, we check that
$$\begin{array}{lll}
a(xt-ty+s)-(xt-ty+s)b&=&x(at-bt)-(at-tb)y+(as-sb)\\
&=&-xv+vy+(as-sb)\\
&=&0.
\end{array}$$ By hypothesis, $xt-ty=-s$, and so we get $$\left(
\begin{array}{cc}
x&s\\
0&y \end{array} \right)\left(
\begin{array}{cc}
0&t\\
0&1 \end{array} \right)=\left(
\begin{array}{cc}
0&ty\\
0&y \end{array} \right)=\left(
\begin{array}{cc}
0&xt+s\\
0&y \end{array} \right)=\left(
\begin{array}{cc}
x&s\\
0&y \end{array} \right)\left(
\begin{array}{cc}
0&t\\
0&1 \end{array} \right).$$ We infer that
$$\left(
\begin{array}{cc}
0&t\\
0&1 \end{array} \right)^2-\left(
\begin{array}{cc}
0&t\\
0&1 \end{array} \right)\in comm^2(A).$$ Furthermore, $A-\left(
\begin{array}{cc}
0&t\\
0&1 \end{array} \right)\in U(T)$. Therefore $A$ is perfectly
clean.

Case $IV$. $a\in J(R), b\in U(S)$, then $A$ is perfectly clean, as
in the preceding discussion.\hfill$\Box$

\vskip4mm A ring $R$ is uniquely weakly bleached provided that for
any $a\in J(R), b\in 1+J(R)$, $l_a-r_b$ and $l_b-r_a$ are both
isomorphisms.

\vskip4mm \hspace{-1.8em} {\bf Corollary 3.6.}\ \ {\it Let $R$ be
local. Then the following are equivalent:} \vspace{-.5mm}
\begin{enumerate}
\item [(1)]{\it $T_2(R)$ is perfectly clean.} \vspace{-.5mm}
\item [(2)]{\it $R$ is uniquely weakly bleached.}
\end{enumerate}
\vspace{-.5mm} {\it Proof.}\ \ It is clear from Theorem 3.5.\hfill$\Box$

\vskip4mm For instance, if $R$ is a commutative local ring or a
local ring with nil Jacobson radical, then $T_2(R)$ is perfectly
clean.

\section{Perfectly $J$-Clean Rings}

\vskip4mm An element $a\in R$ is said to be perfectly $J$-clean
provided that there exists an idempotent $e\in comm^2(a)$ such
that $a-e\in J(R)$. A ring $R$ is perfectly $J$-clean in case
every element in $R$ is perfectly $J$-clean.

\vskip4mm \hspace{-1.8em} {\bf Theorem 4.1.}\ \ {\it Let $R$ be a
ring. Then $R$ is perfectly $J$-clean if and only if}
\vspace{-.5mm}
\begin{enumerate}
\item [(1)]{\it $R$ is quasipolar.} \vspace{-.5mm}
\item [(2)]{\it $R/J(R)$ is Boolean.}
\end{enumerate}
\vspace{-.5mm} {\it Proof.}\ \ Suppose that $R$ is perfectly
clean. Let $a\in R$ is perfectly $J$-clean. Then there exists an
idempotent $e\in comm^2(a)$ such that $w:=a-e\in J(R)$. Hence,
$a-(1-e)=2e-1+w\in U(R)$. Additionally, $(1-e)a=(1-e)w\in
J(R)\subseteq R^{qnil}$. This implies that $a\in R$ is quasipolar.
Furthermore, $a-a^2=(e+w)-(e-w)^2\in J(R)$, and then $R/J(R)$ is
Boolean.

Conversely, assume that $(1)$ and $(2)$ hold. Let $a\in R$. Then
there exists an idempotent $e\in comm^2(a)$ such that $u:=a-e\in
U(R)$. Moreover, $R/J(R)$ is Boolean, and so
$a-a^2=(e+u)-(e+u)^2=u(1-2e-u)\in J(R)$. This shows that
$1-2e-u\in J(R)$; whence, $a-(1-e)=(e+u)-(1-e)=2e-1+u\in J(R)$.
Therefore $R$ is perfectly $J$-clean.\hfill$\Box$

\vskip4mm \hspace{-1.8em} {\bf Corollary 4.2.}\ \ {\it Let $R$ be
a ring. Then the following are equivalent:} \vspace{-.5mm}
\begin{enumerate}
\item [(1)]{\it $R$ is perfectly $J$-clean.} \vspace{-.5mm}
\item [(2)]{\it $R$ is perfectly clean, and $R/J(R)$ is Boolean.}
\vspace{-.5mm}
\item [(3)]{\it $R$ is quasipolar, and $R$ is strongly $J$-clean.}
\end{enumerate}
\vspace{-.5mm} {\it Proof.}\ \ $(1)\Rightarrow (2)$ is obvious
from Theorem 4.1 as every quasipolar ring is perfectly clean.

$(2)\Rightarrow (1)$ For any $a\in R$ there exists an idempotent
$p\in comm^2(a)$ such that $u:=a-p\in U(R)$. As $R/J(R)$ is
Boolean, $\overline{u}=\overline{u}^2$; hence, $u\in 1+J(R)$.
Furthermore, $2\in J(R)$. Accordingly, $a=e+u=(1-p)+(2p-1+u)$ with
$1-p\in comm^2(a)$ and $2p-1+u\in J(R)$, as desired.

$(1)\Rightarrow (3)$ Suppose $R$ is perfectly $J$-clean. Then $R$
is strongly $J$-clean. By the preceding discussion, $R$ is
quasipolar.

$(3)\Rightarrow (1)$ Since $R$ is strongly $J$-clean, $R/J(R)$ is
Boolean. Therefore the proof is complete by the discussion
above.\hfill$\Box$

\vskip4mm \hspace{-1.8em} {\bf Example 4.3.}\ \ {\it Let
$R=T_2({\Bbb Z}_{2^n}) (n\in {\Bbb N})$. Then $T_2(R)$ is
perfectly $J$-clean.}\vskip2mm\hspace{-1.8em} {\it Proof.}\ \ As
$R$ is finite, it is periodic. This shows that $R$ is strongly
$\pi$-regular. Hence, $T_2(R)$ is quasipolar, by [9, Theorem 2.6].
As $J\big({\Bbb Z}_{2^n}\big)=2{\Bbb Z}_{2^n}$, we see that
$R/J(R)\cong {\Bbb Z}_2$ is Boolean. Hence,
$T_2(R)/J\big(T_2(R)\big)$ is Boolean. Therefore the result
follows by Theorem 4.1.\hfill$\Box$

\vskip4mm Recall that a ring $R$ is uniquely strongly clean
provided that for any $a\in R$ there exists a unique idempotent
$e\in comm(a)$ such that $a-e\in U(R)$.

\vskip4mm \hspace{-1.8em} {\bf Proposition 4.4}\ \ {\it Let $R$ be
a ring. Then $R$ is perfectly $J$-clean if and only if}
\vspace{-.5mm}
\begin{enumerate}
\item [(1)]{\it $R$ is perfectly clean.}\vspace{-.5mm}
\item [(2)]{\it $R$ is uniquely strongly clean.}
\end{enumerate}
\vspace{-.5mm} {\it Proof.}\ \ Suppose $R$ is perfectly $J$-clean.
Then $R$ is perfectly clean. Hence, $R$ is strongly clean. Let
$a\in R$. Write $a=e+u=f+v$ with $e=e^2\in comm^2(a),f=f^2\in R,
u\in J(R),v\in U(R), ea=ae$ and $fa=af$. Then $f\in comm(a)$, and
so $ef=fe$. Thus, $e-f=v-u\in U(R)$, and $(e-f)(e+f-1)=0$. This
implies that $f=1-e$, and therefore $R$ is uniquely strongly
clean.

Conversely, assume that $(1)$ and $(2)$ hold. Then $R/J(R)$ is
Boolean. Therefore we complete the proof from Corollary
4.2.\hfill$\Box$

\vskip4mm \hspace{-1.8em} {\bf Corollary 4.5.}\ \ {\it A ring $R$
is uniquely clean if and only if $R$ is abelian perfectly
$J$-clean.}\vskip2mm\hspace{-1.8em} {\it Proof.}\ \ As every
uniquely clean ring is abelian, it is clear from Proposition
4.4.\hfill$\Box$

\vskip4mm \hspace{-1.8em} {\bf Theorem 4.6.}\ \ {\it Let $R$ be a
ring. Then the following are equivalent:} \vspace{-.5mm}
\begin{enumerate}
\item [(1)]{\it $R$ is perfectly $J$-clean.}
\vspace{-.5mm}
\item [(2)]{\it For any $a\in R$, there exists a unique idempotent $e\in comm^2(a)$ such
that $a-e\in J(R)$.}
\end{enumerate}
\vspace{-.5mm} {\it Proof.}\ \ $(1)\Rightarrow (2)$ For any $a\in
R$, there exists an idempotent $e\in comm^2(a)$ such that $a-e\in
J(R)$. Assume that $a-f\in J(R)$ for an idempotent $f\in
comm^2(a)$. Clearly, $e\in comm^2(a)\subseteq comm(a)$. As $f\in
comm(a)$, we see that $ef=fe$. Thus, $(e-f)^3=e-f$, and so
$(e-f)\big(1-(e-f)^2\big)$. But $e-f=(a-f)-(a-e)\in J(R)$. Hence,
$e=f$, as desired.

$(2)\Rightarrow (1)$ is trivial.\hfill$\Box$

\vskip4mm Recall that a ring $R$ is strongly $J$-clean provided
that for any $a\in R$ there exists an idempotent $e\in comm(a)$
such that $a-e\in J(R)$ (cf. [3-4]).

\vskip4mm \hspace{-1.8em} {\bf Corollary 4.7.}\ \ {\it A ring $R$
is perfectly $J$-clean if and only if}\vspace{-.5mm}
\begin{enumerate}
\item [(1)]{\it $R$ is quasipolar;}
\vspace{-.5mm}
\item [(2)]{\it $R$ is strongly $J$-clean.}
\end{enumerate}
\vspace{-.5mm} {\it Proof.}\ \ Suppose $R$ is perfectly $J$-clean.
Then $R$ is strongly $J$-clean. For any $a\in R$, there exists an
idempotent $p\in comm^2(a)$ such that $w:=a-p\in J(R)$. Hence,
$a=(1-p)+(2p-1+w)$ with $1-p\in comm^2(a)$ and $2p-1+w\in U(R)$.
Furthermore, $(1-p)a=(1-p)w\in J(R)\subseteq R^{qnil}$. Therefore,
$R$ is quasipolar.

Conversely, assume that $(1)$ and $(2)$ hold. Since $R$ is
quasipolar, it is perfectly clean. By virtue of [4, Proposition
16.4.15], $R/J(R)$ is Boolean. Therefore the proof is complete by
Corollary 4.2.\hfill$\Box$

\vskip4mm Following Cui and Chen [8], a ring $R$ is called
$J$-quasipolar provided that for any element $a\in R$ there exists
some $e\in comm^2(a)$ such that $a+e\in J(R)$. We further show
that the two concepts coincide with each other. But this is not
the case for a single element. That is,

\vskip4mm \hspace{-1.8em} {\bf Proposition 4.8.}\ \ {\it A ring
$R$ is perfectly $J$-clean if and only if for any element $a\in R$
there exists some $e\in comm^2(a)$ such that $a+e\in
J(R)$.}\vskip2mm\hspace{-1.8em} {\it Proof.}\ \ Let $R$ be
perfectly $J$-clean. Then $R/J(R)$ is Boolean, by Theorem 4.1.
Hence, $\overline{2}^2=\overline{2}$, i.e., $2\in J(R)$. For any
$a\in R$, there exists an idempotent $e\in comm^2(a)$ such that
$a-e\in J(R)$. This implies that $a+e=(a-e)+2e\in J(R)$. The
converse is similar by [8, Corollary 2.3].\hfill$\Box$

\vskip4mm \hspace{-1.8em} {\bf Example 4.9.}\ \ Let $R=\Bbb Z_3$.
Note that $J(R)=0$. Since
$\overline{1}-\overline{1}=\overline{0}\in J(R)$, $\overline{1}$
is perfectly $J$-clean but we can not find an idempotent $e\in R$
such that $\overline{1}+e\in J(R)$, because
$\overline{1}+\overline{0}\notin J(R)$ and
$\overline{1}+\overline{1}=\overline{2}\notin J(R)$.

Further, though $\overline{2}+\overline{1}=\overline{0}\in J(R)$,
we can not find an idempotent $e\in R$ such that
$\overline{2}-e\in U(R)$, because
$\overline{2}-\overline{0}=\overline{2}\notin J(R)$ and
$\overline{2}-\overline{1}=\overline{1}\notin J(R)$.

\vskip4mm \hspace{-1.8em} {\bf Lemma 4.10}\ \ {\it Let $R$ be a
ring. Then $a\in R$ is perfectly $J$-clean if and only if}
\vspace{-.5mm}
\begin{enumerate}
\item [(1)]{\it $a\in R$ is quasipolar.} \vspace{-.5mm}
\item [(2)]{\it $a-a^2\in J(R)$.}
\end{enumerate}
\vspace{-.5mm} {\it Proof.}\ \ Suppose that $a\in R$ is perfectly
$J$-clean. Then there exists an idempotent $e\in comm^2(a)$ such
that $w:=a-e\in J(R)$. Hence, $a-(1-e)=2e-1+w\in U(R)$.
Additionally, $(1-e)a=(1-e)w\in J(R)\subseteq R^{qnil}$. This
implies that $a\in R$ is quasipolar. Furthermore,
$a-a^2=(e+w)-(e-w)^2\in J(R)$.

Conversely, assume that $(1)$ and $(2)$ hold. Then there exists an
idempotent $e\in comm^2(a)$ such that $u:=a-e\in U(R)$. Moreover,
$a-a^2=(e+u)-(e+u)^2=u(1-2e-u)\in J(R)$; hence, $1-2e-u\in J(R)$.
This shows that $a-(1-e)=(e+u)-(1-e)=2e-1+u\in J(R)$. Therefore
$a\in R$ is perfectly $J$-clean.\hfill$\Box$

\vskip4mm \hspace{-1.8em} {\bf Theorem 4.11.}\ \ {\it Let $R$ be a
commutative ring, and let $A\in T_n(R)$. If $2\in J(R)$, then the
following are equivalent:} \vspace{-.5mm}
\begin{enumerate}
\item [(1)]{\it $A\in T_n(R)$ is perfectly $J$-clean.} \vspace{-.5mm}
\item [(2)]{\it Each $A_{ii}\in T_n(R)$ is perfectly $J$-clean.}
\end{enumerate}
\vspace{-.5mm} {\it Proof.}\ \ $(1)\Rightarrow (2)$ is obvious.

$(2)\Rightarrow (1)$ Clearly, the result holds for $n=1$. Suppose
that the result holds for $n-1 (n\geq 2)$. Let $A=\left(
\begin{array}{cc}
a_{11}&\alpha\\
&A_{1}
\end{array}
\right)\in T_n(R)$ where $a_{11}\in R, \alpha\in M_{1\times
(n-1)}(R)$ and $A_1\in T_{n-1}(R)$. Then we have an idempotent
$e_{11}\in R$ such that $w_{11}:=a_{11}-e_{11}\in J(R)$. By
hypothesis, we have an idempotent $E_1\in T_{n-1}(R)$ such that
$W_1:=A_1-E_1\in J\big(T_{n-1}(R)\big)$ and $E_1\in comm^2(A_1)$.
As $2\in J(R)$, $W_1+\big(1-2e_{11}-w_{11}\big)I_{n-1}\in
I_{n-1}+J\big(T_{n-1}(R)\big)\subseteq U\big(T_{n-1}(R)\big)$. Let
$E=\left(
\begin{array}{cc}
e_{11}&\beta\\
&E_1
\end{array}
\right)$, where $\beta=\alpha (E_1-e_{11}I_{n-1})
\big(W_1+(1-2e_{11}-w_{11})I_{n-1}\big)^{-1}$. Then $A-E\in
J\big(T_n(R)\big)$. As $e_{11}\beta+\beta E_1=\beta
(E_1+e_{11}I_{n-1})=\alpha (E_1-e_{11}I_{n-1})(E_1+e_{11}I_{n-1})
\big(W_1+(1-2e_{11}-w_{11})I_{n-1}\big)^{-1}=\beta$, we see that
$E=E^2$.

For any $X=\left(
\begin{array}{cc}
x_{11}&\gamma\\
&X_1
\end{array}
\right)\in comm(A)$, we have $x_{11}\alpha+\gamma
A_1=a_{11}\gamma+\alpha X_1$; hence, $$\alpha
(X_1-x_{11}I_{n-1})=\gamma (A_1-a_{11}I_{n-1}).$$ As $E_1\in
comm^2(A_1)$, we get
$$\begin{array}{ll}
&\gamma(A_1-a_{11}I_{n-1})(E_{1}-e_{11}I_{n-1})\\
=&\alpha (X_1-x_{11}I_{n-1})(E_{1}-e_{11}I_{n-1})\\
=&\alpha (E_{1}-e_{11}I_{n-1})(X_1-x_{11}I_{n-1})\\
=&\beta \big(W_1+(1-2e_{11}-w_{11})I_{n-1}\big)(X_1-x_{11}I_{n-1})\\
=&\beta
(X_1-x_{11}I_{n-1})\big(W_1+(1-2e_{11}-w_{11})I_{n-1}\big).
\end{array}$$
Furthermore,
$$\begin{array}{ll}
&\gamma(A_1-a_{11}I_{n-1})(E_{1}-e_{11}I_{n-1})\\
=&\gamma(E_{1}-e_{11}I_{n-1})\big(E_1+W_1-(e_{11}+w_{11})I_{n-1}\big)\\
=&\gamma(E_{1}-e_{11}I_{n-1})\big(E_1+e_{11}I_{n-1}+(W_1-2e_{11}-w_{11})I_{n-1}\big)\\
=&\gamma\big(E_1-e_{11}I_{n-1}+(E_{1}-e_{11}I_{n-1})(W_1-2e_{11}-w_{11})I_{n-1}\big)\\
=&\gamma
(E_1-e_{11}I_{n-1})\big(W_1+(1-2e_{11}-w_{11})I_{n-1}\big).
\end{array}$$
It follows from $W_1+(1-2e_{11}-w_{11})I_{n-1}\in
U\big(T_{n-1}(R)\big)$ that $\gamma (E_1-e_{11}I_{n-1})=\beta
(X_1-x_{11}I_{n-1})$. Hence, $e_{11}\gamma+\beta
X_{1}=x_{11}\beta+\gamma E_1$, and so $EX=XE$. This implies that
$E\in comm^2(A)$. By induction, $A\in T_n(R)$ is perfectly
$J$-clean for all $n\in {\Bbb N}$.\hfill$\Box$

\vskip4mm \hspace{-1.8em} {\bf Corollary 4.12.}\ \ {\it Let $R$ be
a commutative ring. Then the following are equivalent:}
\vspace{-.5mm}
\begin{enumerate}
\item [(1)]{\it $R$ is strongly $J$-clean.} \vspace{-.5mm}
\item [(2)]{\it $T_n(R)$ is perfectly $J$-clean for all $n\in {\Bbb N}$.}
\item [(3)]{\it $T_n(R)$ is perfectly $J$-clean for some $n\in {\Bbb N}$.}
\end{enumerate}
\vspace{-.5mm} {\it Proof.}\ \ $(1)\Rightarrow (2)$ As $R$ is
strongly $J$-clean, $R/J(R)$ is Boolean. Hence, $2\in J(R)$. For
any $n\in {\Bbb N}$, $T_n(R)$ is perfectly $J$-clean by Theorem
4.11.

$(2)\Rightarrow (3)\Rightarrow (1)$ These are clear from Theorem
4.11.\hfill$\Box$

\vskip4mm Let $R$ be Boolean. As a consequence of Corollary 4.12,
$T_n(R)$ is perfectly $J$-clean for all $n\in {\Bbb N}$.

\vskip4mm \hspace{-1.8em} {\bf Lemma 4.13.}\ \ {\it Let $R$ be a
ring, and let $u\in U(R)$. Then the following are equivalent:}
\vspace{-.5mm}
\begin{enumerate}
\item [(1)]{\it $a\in R$ is perfectly $J$-clean.}
\vspace{-.5mm}
\item [(2)]{\it $uau^{-1}\in R$ is perfectly $J$-clean.}
\end{enumerate}
\vspace{-.5mm} {\it Proof.}\ \ $(1)\Rightarrow (2)$ As in the
proof of Lemma 3.1, $uau^{-1}\in R$ is quasipolar. Furthermore,
$uau^{-1}-\big(uau^{-1}\big)^2=u(a-a^2)u^{-1}\in J(R)$. As in the
proof of Theorem 4.1, $uau^{-1}\in R$ is perfectly $J$-clean.

$(1)\Rightarrow (2)$ is symmetric.\hfill$\Box$

\vskip4mm We end this paper by showing that strong $J$-cleanness
of $2\times 2$ matrix ring over a commutative local ring can be
enhanced to perfect $J$-cleanness.

\vskip4mm \hspace{-1.8em} {\bf Theorem 4.14.}\ \ {\it Let $R$ be a
commutative local ring, and let $A\in M_2(R)$. Then the following
are equivalent:} \vspace{-.5mm}
\begin{enumerate}
\item [(1)]{\it $A$ is perfectly $J$-clean.} \vspace{-.5mm}
\item [(2)]{\it $A$ is strongly $J$-clean.} \vspace{-.5mm}
\item [(3)]{\it $A\in J\big(M_2(R)\big)$, or $I_2-A\in J\big(M_2(R)\big)$, or the equation $x^2-tr(A)x+det(A)=0$ has a root in $J(R)$ and a root in $1+J(R)$.}
\end{enumerate}
\vspace{-.5mm} {\it Proof.}\ \ $(1)\Rightarrow (2)$ is trivial.

$(2)\Rightarrow (3)$ is proved by [4, Theorem 16.4.31].

$(3)\Rightarrow (1)$ If $A\in J\big(M_2(R)\big)$ or $I_2-A\in
J\big(M_2(R)\big)$, then $A$ is perfectly $J$-clean. Otherwise, it
follows from [4, Theorem 16.4.31 and Proposition 16.4.26] that
there exists a $U\in GL_2(R)$ such that
$$UAU^{-1}=\left(
\begin{array}{cc}
\alpha&\\
&\beta \end{array} \right)=\left(
\begin{array}{cc}
0&\\
&1 \end{array} \right)+\left(
\begin{array}{cc}
\alpha&\\
&\beta-1 \end{array} \right),$$ where $\alpha\in J(R),\beta\in
1+J(R)$. For any $X\in comm(UAU^{-1})$, we have $X\left(
\begin{array}{cc}
\alpha&\\
&\beta \end{array} \right)=\left(
\begin{array}{cc}
\alpha&\\
&\beta \end{array} \right)X$; hence, $\beta X_{12}=\alpha X_{12}$.
This implies that $X_{12}=0$. Likewise, $X_{21}=0$. Thus,
$$X\left(
\begin{array}{cc}
0&\\
&1 \end{array} \right)=\left(
\begin{array}{cc}
0&\\
&1 \end{array} \right)X,$$ and so $\left(
\begin{array}{cc}
0&\\
&1 \end{array} \right)\in comm^2(UAU^{-1})$. As a result,
$UAU^{-1}$ is perfectly $J$-clean, and then so is $A$ by Lemma
4.13.\hfill$\Box$

\vskip4mm \hspace{-1.8em} {\bf Corollary 4.15.}\ \ {\it Let $R$ be
a commutative local ring. Then the following are equivalent:}
\vspace{-.5mm}
\begin{enumerate}
\item [(1)]{\it $M_2(R)$ is perfectly clean.} \vspace{-.5mm}
\item [(2)]{\it For any $A\in M_2(R)$, $A\in GL_2(R)$, or $I_2-A\in GL_2(R)$, or $A\in M_2(R)$ is perfectly $J$-clean.}
\end{enumerate}
\vspace{-.5mm} {\it Proof.}\ \ $(1)\Rightarrow (2)$ is proved by
Theorem 3.2, [4, Corollary 16.4.33] and Theorem 4.14.

$(2)\Rightarrow (1)$ is obvious.

\vskip10mm \no {\Large\bf Acknowledgements} \vskip4mm \no This
research was supported by the Natural Science Foundation of
Zhejiang Province (LY13A0 10019) and the Scientific and
Technological Research Council of Turkey (2221 Visiting Scientists
Fellowship Programme).

\end{document}